\documentclass[10pt]{article}

\usepackage[a4paper, total={6in, 9in}]{geometry}
\usepackage{color}
\usepackage{graphicx}
\usepackage{latexsym}
\usepackage{amsfonts,amsmath,amssymb}

\newcommand{\qed}{\mbox{$\Diamond$}\vspace{\baselineskip}}

\begin{document}

\author{Istv\'an Mez\H{o}\thanks{The research of Istv\'an Mez\H{o} was supported by the Scientific Research Foundation of Nanjing University of Information Science \& Technology, and The Startup Foundation for Introducing Talent of NUIST. Project no.: S8113062001}\\Department of Mathematics\\Nanjing University of Information Science and Technology\\
Nanjing, 210044, P. R. China}

\title{Asymptotics of the modes of the ordered Stirling numbers}
\maketitle

\begin{abstract}It is known that the $S(n,k)$ Stirling numbers as well as the ordered Stirling numbers $k!S(n,k)$ form log-concave sequences. Although in the first case there are many estimations about the mode, for the ordered Stirling numbers such estimations are not known. In this short note we study this problem and some of its generalizations.
\end{abstract}

\section{Introduction}

It is a classical result that for a fixed $n$ the $S(n,k)$ Stirling numbers form strictly log-concave sequences (SLC), that is, for any $n>2$ it holds true that
\[S(n,1)<S(n,2)<\cdots<S(n,K_n^*)\ge S(n,K_n^*+1)>S(n,K_n^*+2)>\cdots>S(n,n)\]
for some index $K_n^*$. It is conjectured by Wegner that $S(n,K_n^*)> S(n,K_n^*+1)$, that is, the maximum is unique for $n>2$ (see also \cite{CP,KMS}). There are good estimations for the mode $K_n^*$, it asymptotically equals to $n/\log(n)$ (see \cite{MC} for a short proof).

It is somewhat less known that the sequence $k!S(n,k)$ is also log-concave \cite{MezoThesis}. Let us denote its mode (or the smallest if there are two) by $M_n$. We are going to prove that as $n\to\infty$
\begin{equation}
M_n\sim\frac{n}{2\log2}.\label{Masy}
\end{equation}

We prove the corresponding results for the $r$-Stirling numbers and Whitney numbers -- both of them are generalizations of the classical Stirling numbers. Thus, specializing these, we will have two independent proof of \eqref{Masy}.

\section{The $r$-Stirling numbers}

Darroch \cite{Darroch} proved that if for a real sequence $a_0,a_1,\dots,a_n$ the attached polynomial $p(x)=a_0+a_1x+\cdots+a_nx^n$ has only real zeros and $p(1)>0$, then the sequence $a_0,a_1,\dots,a_n$ is SLC, and the (smallest) mode $M$ can be located as
\[\left|M-\frac{p'(1)}{p(1)}\right|<1.\]

We apply this result first to the $r$-Stirling numbers. The $S_r(n,k)$ Stirling numbers \cite{Broder} count the partitions of an $n+r$ element set into $k+r$ subsets such that $r$ distinguished elements are restricted to be in different subsets. If $r=0$ we get back the classical Stirling numbers. It was proven by the present author that $(k+r)!S_r(n,k)$ is SLC \cite{MezoThesis}, but the asymptotic behavior of the mode -- what we will denote by $M_{n,r}$ -- was not studied. Here we prove that $M_{n,r}$ asymptotically equals to $M_n$:
\begin{equation}
M_{n,r}\sim\frac{n}{2\log 2}.\label{rStirmode}
\end{equation}

The $r$-Fubini polynomials
\[F_{n,r}(x)=\sum_{k=1}^n(k+r)!S_r(n,k)x^k\]
were studied in detail \cite{MezoThesis}. It was proven that $F_{n,r}(x)$ has only real zeros (and they lie in the interval $]-1,0]$). So the sequence $((k+r)!S_r(n,k))_{k=1}^n$ is SLC. Since $F_{n,r}(x)$ satisfies the recurrence
\[F_{n,r}(x)=x[(r+1)F_{n-1,r}(x)+(1+x)F_{n-1,r}'(x)]+rF_{n-1,r}(x),\]
it can easily be seen that $F_{n,r}'(x)=\frac{F_{n+1,r}(x)-(2r+1)F_{n,r}(x)}{2}$, so Darroch's theorem says that
\begin{equation}
\left|M_{n,r}-\frac{F_{n+1,r}-(2r+1)F_{n,r}}{2F_{n,r}}\right|<1.\label{Mnrest}
\end{equation}
Here $F_{n,r}=F_{n,r}(1)$. The asymptotics of the $F_{n,r}$ numbers was not examined in \cite{MezoThesis}, so we do it here. The exponential generating function of these numbers is
\[\sum_{n=0}^\infty F_{n,r}\frac{t^n}{n!}=\frac{r!e^{rt}}{(2-e^t)^{r+1}}.\]
The smallest pole of this function is at $t=\log 2$ with order $r+1$. Analyzing carefully the Laurent series around this pole and applying the standard methods of the saddle point method \cite{Wilf} we can show by a lengthier computation that
\begin{equation}
F_{n,r}\sim\frac{1}{2\log^{r+1}(2)}\frac{n!n^r}{\log^n(2)}.\label{Fasymp}
\end{equation}
(The details are left to the reader.) This is a generalized version of a theorem of Vellemand and Call who proved the particular case
\[F_{n,0}\sim\frac{1}{2\log2}\frac{n!}{\log^n(2)}.\]

Substituting \eqref{Fasymp} into \eqref{Mnrest} the asymptotic behavior \eqref{rStirmode} already follows. Taking $r=0$ we get back \eqref{Masy}, too.

\section{The Whitney numbers}

The $W_m(n,k)$ Whitney numbers are another generalizations of the Stirlings via lattice theory \cite{Benoum,Benoum2} such that $W_1(n,k)=S(n+1,k+1)$. These sequences are SLC \cite{Benoum3} as well as $k!W_m(n,k)$ \cite[Theorem 6]{Benoum2}. The asymptotics of the mode of the former sequence was determined by Benoumhani \cite{Benoum}; it is $\frac{n}{m\log n}$. Here we deal with the mode of $k!W_m(n,k)$. The polynomials
\[F_m(n,x)=\sum_{k=1}^nk!W_m(n,k)x^k\]
satisfy the recursion \cite{Benoum2}
\[F_m(n,x)=(x+1)F_m(n-1,x)+(x^2+mx)F_m'(n-1,x),\]
and the exponential generating function of $F_m(n)=F_m(n,1)$ is
\[\sum_{n=0}^\infty F_m(n)\frac{t^n}{n!}=\frac{e^t}{1-\frac1m(e^{mt}-1)}.\]
This function has a simple pole at $t=\log(m+1)/m$. Analyzing the Laurent expansion around this point and perform the standard steps of the saddle point method after some algebra we get the following asymptotic behavior of the $F_m(n)$ numbers:
\[F_m(n)=\frac{m(1+m)^{\frac1m-1}}{\log(m+1)}\frac{m^nn!}{\log^n(m+1)}.\]
From this asymptotics it already follows that, denoting the mode of $k!W_m(n,k)$ by $W_{n,m}$,
\[W_{n,m}\sim\frac{m}{(m+1)\log(m+1)}n.\]
Setting $m=1$ we have \eqref{Masy}, again.

\end{document}